\documentclass[12pt]{article}
\usepackage{newtxtext}          %
 \usepackage{newtxmath}
 \usepackage[T1]{fontenc}
 \usepackage[utf8]{inputenc}
\usepackage[numbers,sort&compress]{natbib}
\usepackage{amsmath}
\usepackage{hyperref}
\usepackage{graphicx}
\usepackage{geometry}

\hypersetup{
    pdfnewwindow=true,          
    colorlinks=true,            
    linkcolor=blue,            
    citecolor=blue,             
    filecolor=magenta,          
    urlcolor=cyan,              
    pdfauthor=D\"untsch \& Dzik,          
}

 \parskip=\medskipamount
\author{Viktor Marek\footnote{University of Kentucky, Lexington,  \url{marek@cs.uky.edu }},\; Ewa Or{\l}owska\footnote{National Institute of Telecommunications, Warsaw, \url{stella.ewa.orlowska@gmail.com}}~ \& Ivo D{\"u}ntsch\footnote{Brock University, St. Catharines, \url{duentsch@brocku.ca}}}
\title{The Evolution of Rough Sets \\ 1970s -- 1981}
\date{}

\begin{document}

\maketitle

\begin{abstract}
\noindent In this note research and publications by Zdzis\l{}aw Pawlak and his collaborators from 1970s and 1981 are recalled. Focus is placed on the sources of inspiration which one can identify on the basis of those publications. Finally, developments from 1981 related to rough sets and information systems are outlined.
\end{abstract}

\section{Early history}

The original work in Poland on computers was made mostly by electrical engineers and mathematicians. In this context, in the early 1960s a group of researchers and students formed a seminar devoted to computer science, formal linguistics, and logic. The seminar met on a weekly basis at the Mathematical Institute of the Polish Academy of Sciences (PAS). It was devoted to a wide spectrum of problems. The common ``Leitmotiv'' was the use of computers as tools to represent and process problems stemming of formal linguistics, mathematical logic, automated theorem proving, models of computation, and other areas on the borderline of computer science and related theories.

The seminar was organized by Robert Bartoszy\'{n}ski (1933–1998), a statistician, Andrzej Ehrenfeucht, a logician,  and Zdzis\l{}aw Pawlak (1926-2006), a computer engineer. Irena Bellert (1919-2017), a linguist,  was an active participant of the seminar. The first and second author of this note took part in the seminar and contributed their presentations. 

The effect of the work by this group resulted in the creation of a nucleus of the team that subsequently worked at Warsaw University, mostly within the section of the Mathematical Institute. The head of this group was Professor Helena Rasiowa (1917-1994), a well-known logician who specialized in algebraic methods applied to logic. The group of individuals involved in this group over the years included A. Jankowski, B. Konikowska, G. Mirkowska-Salwicka, A. Salwicki, and A. Skowron, among others. Individuals  from outside the University included W. Lipski and M. Jaegermann.

However, in the early 1970s it became clear that research in the area of information systems had run its course. Several variations had been invented and subsequently investigated, papers written, even implementations proposed. But it was a growing consensus that new ideas were needed --  within the same framework, but going beyond the current proposals.

\section{Inspiration from Computer Science}

The concept of fuzzy sets was proposed by Lotfi \citet{zadeh65,zadeh75} in 1965 for dealing with uncertainty of information resulting from its incompleteness.\footnote{Independently, a similar concept was proposed by \citet{klaua65,klaua66}.} A set is \emph{fuzzy} whenever the range of its membership function is the real interval [0,1]. Fuzzy sets provide a means for quantitative representation of uncertain information.

On the way for defining a formal model of collection of data, the relational database model proposed by E. F. \citet{cod70} was an important step. In this
model, a database is an $n$-ary relation, say $R$, represented as an array of rows such that each row is an $n$-tuple of $R$, the ordering of rows is immaterial and all rows are distinct. However, by treating a row as an object and elements in the row as its properties, it is natural to assume that not necessarily only a single object has the properties listed in the row. In such a case, these objects cannot be differentiated with the properties in question. As a consequence, a precise definition of the set of those objects in terms of their  properties cannot be given.

\section{Inspiration from Logic and Linguistics}

Among the important formal tools for coping with imprecision are multiple-valued logics originated by the Polish logician and philosopher Jan \L{}ukasiewicz (1879-1956). In 1920 \citet{luk1920} proposed a three-valued logic with an intermediate value between true and false.\footnote{For an English translation see
\cite{luk70}.} The logic was generalized to $n$-valued logic for finite $n$, and to infinitely-many-valued $\aleph_{0}$-valued logic by \citet{lt30} in 1930.

A similar idea expressed in algebraic form can be found in the PhD dissertation of Emil Post (1879-1954), a shortened version of which was published in the American Journal of Mathematics \cite{post1921}.

\section{Inspiration from Mathematics}

In the 1970s Rudolf  Wille (1937-2017) -- a professor of General Algebra at the Technical University of  Darmstadt -- visited several times algebraists at the Technical University of Warsaw and gave a series of lectures on the application of lattice theory to organization of collections of data, based on objects and their properties/attributes. His ideas were presented in \cite{wille82}. This paper gave rise to Formal Concept Analysis (FCA). In FCA the model of data is referred to as a formal context: A \emph{formal context} is a triple (objects, properties, information relation), where an information relation consists of pairs of the form (object, property) with the intuition that the object in the pair possesses the property. Subsets of objects are interpreted as extents of concepts, and subsets of properties are interpreted as intents of concepts. In logic, the extent of a concept is a set of objects to which the concept applies, and the intent of a concept is a set of properties or qualities associated with the concept.

Given a formal context $(X, Y, I)$ with a set $X$ of objects, a set $Y$ of properties of objects, and an information relation $I \subseteq X \times Y$, a \emph{formal concept} of the context is a pair $(A, B)$ such that $A$ is a subset of $X$, $B$ is a subset of $Y$, and $A$ is included in $[[I]](B)$ if and only if
$B$ included in $[[I^{-1}]](A)$, where
\begin{align*}
[[I]](B) &= \{x \in X: (\forall y \in Y)[y \in B \text{ implies } x \mathrel{I}y\} \\
[[I^{-1}]](A) &= \{y \in Y: (\forall x \in X)[x \in A \text{ implies } y\mathrel{I}x\}.
\end{align*}
From a formal context one can derive a formal concept lattice \cite{wille82}.

\section{Inspiration from Philosophy}

Marian Przelecki (1923-2013), a Polish professor of  Philosophical Logic, drew the second author's attention to the problem of vagueness of
concepts and Leibniz's principle of indiscernibility of identicals: Leibniz's principle states that no two distinct objects can have exactly the same properties. The annual conferences on the history of logic held at the Jagiellonian University in Krak\'{o}w gave us an opportunity to discuss these issues which later enabled defining the indiscernibility relation and a lower (resp. upper) approximation of  sets in an  information system. A concept is vague whenever its extent is intrinsically uncertain which means that there are objects whose membership in the concept is intrinsically uncertain. The classical example of such a concept is a heap of grains. Early studies of the problem can be found in \citet{russell23} and \citet{black37,black63}. 

\section{Description language and descriptive systems}

Inspiration from the sources listed above led to the concept of descriptive language and a descriptive system presented in \cite{paw73}. Atomic expressions of the descriptive language, referred to as \emph{elementary descriptors}, were of the form \emph{(attribute,value)}. Compound expressions were formed from atomic expressions with Boolean propositional connectives. A \emph{descriptive system} is a triple \emph{(Objects, Elementary Descriptors, Description Relation)}, where  a description relation is a set of pairs of the form \emph{(object, elementary descriptor)}. It was a predecessor of the concept of information system.

In 1974 Zdzis\l{}aw Pawlak, together with his collaborators, proposed a formal model of collections of data referred to as an information storage and retrieval system. A discussion and some properties of these systems can be found in \citet{jae75}, \citet{lipski74}, and \citet{mp74}. The simpler name ``information system'' was adopted from the journal ``Information Systems'', published by Pergamon Press since 1975. In the late 1970s the first author of this note spent time out of Poland with limited contact to Zdzis\l{}aw Pawlak, but he returned in March 1980 just in time to contribute to research on rough sets in 1981.

\section{1981 - The birth of rough sets}

In 1981, Zdzis\l{}aw Pawlak invited to Warsaw Erhard  Konrad (PhD) of the Technical University of Berlin ``to start research in a new field'' -- as he said. The results of the joint work were reported in the two papers by \citet*{kop81a,kop81b} as well as in the article \cite{kop81c} which was intended for submission to the 7$^{th}$ International Joint Conference on Artificial Intelligence. However, owing to the financial difficulties to obtain funds  -- Polish currency was not convertible at that time -- it was never submitted. In the second of those papers the lower (resp. upper) description of a set of objects was defined as an appropriate expression of the description language and the set providing semantics of that expression was named the lower (resp. upper) approximation of the given set.

At that stage, the only missing step on the way to the concept of a rough set was to establish a framework for operating with sets and their approximations. The idea was very natural: Assign to a given set of attributes both interior and closure operators on the set of rows (also referred to as records). In \cite{paw81a,paw82} the interior (resp. closure) operator was called lower (resp. upper) approximation. These operators measured, in a sense, the similarity of records in the database to a given set of records. This was stated precisely in the ICS PAS report \cite{paw81a} published in 1981:  First, an \emph{approximation space} is defined as a pair $(U, R)$ where $U$ is a set of objects and $R$ is an equivalence relation on $ $U. Its equivalence classes are referred to as \emph{elementary sets}. Then, given a subset $X$ of $U$, its \emph{lower} (resp. \emph{upper}) \emph{approximation} was defined as the union of elementary sets included in $X$ (resp. the union of elementary sets intersecting $X$). Next, the relations of approximate inclusions and approximate equalities of sets were presented in terms of their approximations. No explicit definition of rough set was given.

\begin{figure}[htb]
 \caption{A rough set \cite[p. 343]{paw82}}\label{fig:rough}
  \centering
  \includegraphics[width=0.5\textwidth]{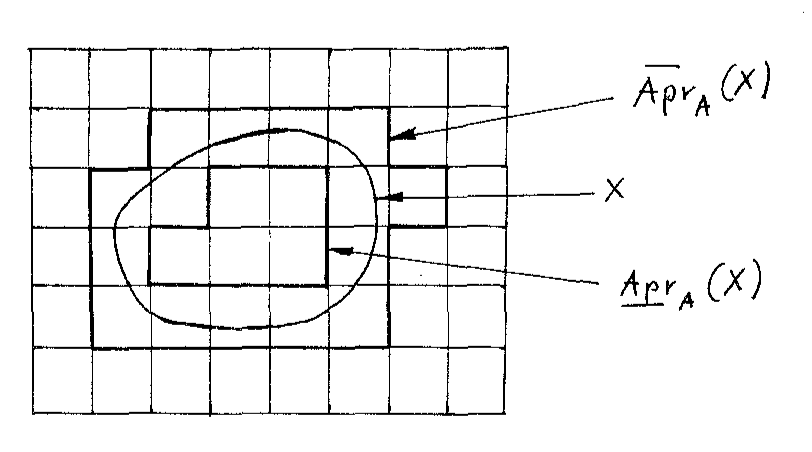}
 \end{figure}

Finally in \cite{paw82} (submitted in 1981), following Marian Przelecki, the relation $A$ in an approximation space $(U, A)$ was named \emph{indiscernibility relation}, and the concept of rough set was presented in the form of a drawing showing, given the partition of set U determined by elementary sets and a subset X of U, its upper and lower approximations provided by the partition. Figure \ref{fig:rough} shows the original drawing from \cite{paw82}.

Altogether, in 1981 several papers by Zdzis\l{}aw Pawlak were published on information systems, rough sets and related topics. The titles of  ICS PAS Reports published by Zdzis\l{}aw Pawlak at that time are: 
\begin{itemize}
\item Classification of objects by means of attributes \cite[ICS PAS Report 429]{paw81c},
\item Rough sets. Basic notions \cite[ICS PAS Report 431]{paw81a},
\item Rough relations \cite[ICS PAS Report 435]{paw81d},\footnote{Also appeared in \cite{paw_roughrel}.}
\item About conflicts \cite[ICS PAS Report 451]{paw81f}
\item Rough functions \cite[ICS PAS Report 467]{paw81e},
\end{itemize}
Furthermore, \citeauthor{il81} published the ICS PAS Report 446  ``The relational model and cylindric algebras'' \cite{il81}. In December 1981 an ICS PAS Report (427) \cite{ICS427} was published of which Part 1 lists a bibliography on the theory of information systems initiated by Zdzis\l{}aw Pawlak and the first author of the present note.

From the logic point of view the modal logic \textbf{S5} is an adequate abstract framework which provides both Hilbert style axioms and  frame semantics (also referred to as possible world semantics) for the class of rough sets.

In the following years a number of papers on applications of rough sets and generalizations of the concept of rough set were published. It is now a thriving area of computer science.

\section*{Acknowledgment}

This note was conceived and drafted by the first two authors and edited by the third author.

\section*{References}
\renewcommand*{\refname}{}

\end{document}